\documentclass[leqno]{article}% 11pt,
\usepackage{amsmath,amstext, amsthm, amssymb}
\numberwithin{equation}{section}
\newtheorem{thm}{Theorem}[section]
\newtheorem{lem}[thm]{Lemma}

\newtheorem{cor}[thm]{Corollary}

\newcommand{\be}{\begin{equation}}
\newcommand{\ee}{\end{equation}}
\newcommand{\ba}{\begin{array}}
\newcommand{\ea}{\end{array}}

\newcommand{\al}{\alpha}

\renewcommand{\th}{\theta}

\newcommand{\la}{\lambda}

\newcommand{\fibon}[2]{\genfrac{(}{)}{0pt}{}{#1}{#2}_{\mathbb{F}}}
\newcommand{\luc}[2]{\genfrac{(}{)}{0pt}{}{#1}{#2}_{\mathbb{L}}}

\newcommand{\bea}{\begin{eqnarray}}
\newcommand{\eea}{\end{eqnarray}}
\newcommand{\Sum}{\sum_{n=0}^\infty}

\title{One Parameter Generalizations of the Fibonacci and Lucas Numbers}
    \author{Mourad E.H. Ismail
\\ Department of Mathematics \\ University of Central Florida
\\  Orlando, FL 32816\\
email address:  ismail@math.ucf.edu}
\begin{document}
\maketitle
\begin{abstract}
We give one parameter generalizations of the Fibonacci and Lucas
numbers denoted by $\{F_n(\th)\}$ and $\{L_n(\th)\}$,
respectively. We evaluate the Hankel determinants with entries
$\{1/F_{j+k+1}(\th): 0\le i,j \le n\}$ and $\{1/L_{j+k+1}(\th):
0\le i,j\le n\}$. We also find the entries in the inverse of
$\{1/F_{j+k+1}(\th): 0\le i,j \le n\}$ and show that all its
entries are integers. Some of the identities satisfied by the
Fibonacci and Lucas numbers are extended to more general numbers.
All integer solutions to three diophantine equtions related to the
Pell equation are also found.
\end{abstract}
\vskip4pt
\noindent
 \noindent{\bf Running Title}.  Generalized Fibonacci and Lucas Numbers
 \vskip2pt
 \noindent

\bigskip

% \noindent {\bf Reference}:

\bigskip
\section{Introduction}
The Hilbert matrix $H_n$ has entries $1/(i+j+1): 0\le i,j \le n$.
It is well-known that, for all $n$,  $H_n$ is non singular and the
elements of its inverse matrix are all integers, see for example
\cite{Cho}. The determinant of $H_n$ has a closed form expression
which shows that the determinant is very small for large $n$. This
is important in numerical analysis because the smaller the
determinant, the larger the condition number becomes and computing
the inverse numerically becomes unstable. Many other applications
of the Hilbert matrix are in \cite{Cho}.

The Fibonacci numbers have many interesting properties and appear
in many areas of mathematics. \cite{Kos}, \cite{Vaj}. One
unexpected result is due to Richardson who showed in \cite{Ric}
that the ``Filbert matrix" is also non singular and its inverse
has only integer entries. The $i,j$ entry of the Filbert matrix is
$1/F_{i+j+1}$ where $0 \le i,j \le n$ and $\{F_n:n \ge 1\}$ are
the Fibonacci numbers.

One way to compute the determinant  of $H_n$ is to note that
it is the Hankel determinant associated with a constant weight
function supported on $[0, 1]$.
Berg \cite{Ber} observed that the reciprocals of the Fibonacci
numbers form a moment sequence of a special little $q$-Jacobi weight,
\cite[\S 18.4]{Ismbook}. He used Lemma 1.1, to be stated
below,  to   prove Richardson's result.

Recall that the Chebyshev polynomials of the first and second
kinds are
\bea
T_n(\cos \th) = \cos (n\th), \quad
U_n(\cos \th) = \frac{\sin ((n+1)\th)}{\sin \th},
 \eea
respectively. Askey \cite{Ask1}, \cite{Ask2} observed that the Fibinacci numbers
$\{F_n\}$ and the
Lucas numbers $\{L_n\}$ are related to the Chebyshev polynomials via
\bea
F_{n+1} = (-i)^{n} U_{n}(i \sinh (\th_0)), \quad
L_{n+1} = 2 (-i)^n T_{n}(i\sinh \th_0))
\eea
where $\th_0 >0$ and $\sinh \th_0 = 1/2$.

The purpose of this paper is to give one parameter generalizations
of the Fibonacci and Lucas numbers. Our generalization  comes from
the Chebyshev polynomials of the first and second kinds.  Our
generalizations of the  Fibonacci and Lucas numbers satisfy the
recurrence relation %
\bea \label{eqfibrec} y_{n+1}(\th) = 2\sinh
\th y_n(\th) + y_{n-1}(\th) %
\eea %
and the initial conditions \eqref{eqincondfn} and
\eqref{eqincondln}. We generalize Richardson's result by replacing
the Fibonacci numbers by our generalized Fibonacci numbers. This
will be done in \S2. In \S 3 we introduce the generalized Lucas
numbers and study some of their properties. We also give a closed
form evaluation of  a Hankel determinant whose elements are
reciprocals of Lucas numbers.

The book \cite{Kos} contains many results on Fibonacci and Lucas
numbers with detailed proofs. On the other hand we found Vajda's
book \cite{Vaj} to be very comprehensive but concise.  In \S 4 we
extend some of the properties of the Fibonacci and Lucas numbers
to our numbers. We have only included a sample of the identities
involving $\{F_n(\th)\}$ and $\{L_n(\th)\}$. There are many other
relationships involving the Fibonacci and Lucas numbers which
extend to our more general sequences $\{F_n(\th)\}$ and
$\{L_n(\th)\}$ but we made no attempt to include them. In \S 5 we
 describe all integer solutions to
$$ y^2 -kxy -x^2 = \pm 1, $$
for a given integer $k > 1$.  We also characterize all integers
$n$ for which $n^2 (1+k^2) \pm 4$ is a perfect square when $k$ is
odd. When $k=1$ these results reduce to known facts involving the
Fibonacci numbers.

The connection between Fibonacci numbers, hyperbolic functions,
and Chebyshev polynomials was observed but some how never fully
exploited, see for example \cite[Chapter 11]{Vaj}, and
\cite{Mel:Sha}. Another recently development is due 
to Kalman and Mena \cite{Kal:Men} who treated 
sequences which satisfy the three term recurrence relation 
$$
y_{n+1}= a y_n + b y_{n-1},
$$ 
under general initial conditions. They derived many of the properties 
that their generalized 
sequence share with the Fibonacci or Lucas numbers.  Our numbers  
being less general than the Kalman-Mena numbers have additional 
properties. For example the inverse matrix to $1/L_{1+i+j}$ does not 
have integer coefficients, while the inverse matrix to $1/F_{i+j+1}$ as 
well as $1/F_{i+j+1}(\theta)$ have integer entries, $\{F_n(\th)$ 
being our generalization of the Fibonacci number. 
Some of the other refined properties involving congruences and 
integer points on algebraic curves or surfaces do not extend  to the 
very general setting of Kalman and Mena. 

 We now come to Lemma 1.1.

Let $\mu$ be a  measure whose moments, $s_0\ne 0$, $s_n
:=\int_{\mathbb{R}} x^n d\mu(x)$ exist for all $n=0,1,\dots$, and
let $\{p_n(x)\}$ be the sequence of polynomials orthogonal with
respect to $\mu$, that is 
\bea 
\label{eqorth} 
\int_{\mathbb{R}}
p_m(x)p_n(x) d\mu(x) = \zeta_n \delta_{m,n}, \quad \zeta_n \ne 0
\eea for $n =0, 1, 2, \dots$. We shall always normalize $\mu$ by
$\zeta_0 = 1$, so that $\mu$ has a unit total mass.
 The corresponding Hankel matrix and Hankel
determinant are  
\bea 
\qquad \;\; H_n =\left(\ba{cccc}
s_0 & s_1 & \dots & s_n\\
s_1 & s_2& \dots & s_{n+1} \\
 \vdots &\vdots  & \dots& \vdots\\
s_n & s_{n+1} & \dots& s_{2n}
\ea
\right), \; D_n =  \left|\ba{cccc}
s_0 & s_1 & \dots & s_n\\
s_1 & s_2& \dots & s_{n+1} \\
 \vdots &\vdots  & \dots& \vdots\\
s_n & s_{n+1} & \dots& s_{2n}
\ea
\right|,
\eea
respectively, and  n = 0, 1, \dots.
The kernel polynomials are
\bea
K_n(x,y) = \sum_{k=0}^n p_k(x) p_k(y)/\zeta_n.
\eea
\begin{lem}
Let
\bea
K_n(x,y) = \sum_{j, k=0}^n a_{j,k}(n) \, x^j \, y^k.
\eea
Then $a_{j,k}(n) = a_{k,j}(n)$ and the matrix $A_n$ whose entries
are $\{a_{j,k}(n)\}$ is the inverse of $H_n$.
\end{lem}
Lemma 1.1 is in the paper \cite{Tra:Wid} by Tracy and Widom and in
Berg's paper \cite{Ber}.

\section{Generalized Fibonacci Numbers}

Consider the Chebyshev polynomials of the second kind %
\bea
\begin{gathered}
U_n (i\sinh \th) = U_n(\cos(\pi/2-i\th) =
\frac{e^{i(\pi/2-i\th)(n+1)} - e^{-i(\pi/2-i\th)(n+1)}}
{e^{i(\pi/2-i\th)} - e^{-i(\pi/2-i\th)}}\\
= i^n \; \frac{e^{(n+1)\th} + (-1)^ne^{-(n+1)\th}} {e^\theta +
e^{-\th}}.
\end{gathered}
% \nonumber
\eea
Set
\bea
\label{eqdefnFn}
F_{n+1}(\th) = (-i)^n\;  U_n(\sinh(i\th))
 = \frac{e^{(n+1)\th} + (-1)^ne^{-(n+1)\th}}{e^\theta + e^{-\th}}.
\eea %
The explicit representation of the Chebyshev polynomials of the
second kind, see for example \cite[\S 4.5]{Ismbook} leads to %
\bea %
\label{eqfnascheby} %
F_n(\th) = \sum_{k=0}^{\lfloor{n/2}\rfloor} \binom{n+1}{2k+1}
\sinh^{n-2k}( \th) \cosh^{2k}(\th).%
 \eea

Choose $\th_0 >0$ so that $\cosh \th_0 = \sqrt{5}/2$. Thus $\sinh
\th_0 =1/2$ and $e^{\th_0}= \phi$ in Berg's notation in
\cite{Ber}. Clearly $e^{-\th_0} = (\sqrt{5}-1)/2 - \hat{\phi}$ in
Berg's notation. Thus $F_n(\th_0)= F_n$, $n =1, 2, \dots$, the
Fibonacci
sequence. Moreover %
\bea F_n(\th) = e^{n\th}\; \frac{1- (-1)^n
e^{-2n\th}}{e^\th + e^{-\th}}. %
\eea %
For positive integer $\al$ we
have %
\bea \frac{F_\al(\th)}{F_{n+\al}(\th)} = e^{-n\th} \;
\frac{1-(-e^{-2\th})^\al}{1- (-e^{-2\th})^{n+\al} }.
\eea
With
\bea \label{eqdfnq} q= - e^{-2\th} %
\eea %
we arrive at %
\bea
\label{eqfnasq} F_n(\th) = e^{(n-1)\th}\; \frac{1-q^n}{1-q}.%
 \eea
Formula \eqref{eqfnasq} enables us to extend the definition of
$F_n(\th)$ to nonpositive  values of $n$. This agrees with
defining $F_n(\th)$ for $n \le 0$ from \eqref{eqfibrec} and the
initial conditions \eqref{eqincondfn} below. Indeed it is easy to
see that %
\bea F_{-n}(\th) = (-1)^{n-1}F_n(\th). \label{eqf-n}
\eea
 {} From \eqref{eqfnasq} it follows that
 \bea \label{eqfib/fib} \frac{F_\al(\th)}{F_{n+\al}(\th)}
= (1-q^\al)\sum_{k=0}^\infty (q^k/e^\th)^n \, q^{\al k}.
\eea

Now use $\fibon{n}{k}$ to denote the binomial coefficient relative to
$\{F_n(\th)\}$, that is
\bea
\label{eqdeffibbinom}
\fibon{n}{0} := 1, \quad \fibon{n}{k}=
\frac{F_n(\th) F_{n-1}(\th)\dots F_{n-k+1}(\th)}
{F_1(\th)\; F_2(\th)\;  \dots \; F_k(\th)}.
\eea
\begin{thm}
We have
\bea
\label{eqrecbinom}
\fibon{n}{k} = F_{k-1}(\th) \fibon{n-1}{k}
+ F_{n-k+1}(\th) \fibon{n-1}{k-1}.
\eea
\end{thm}
\begin{proof}
It is easy to write the right-hand side of \eqref{eqrecbinom}
in the form
\bea
\begin{gathered}
\frac{1}{F_n(\th)}\fibon{n}{k}[F_{k-1}(\th)F_{n-k}(\th)
+ F_k(\th) F_{n-k+1}(\th)] \\
= \frac{(1-q)^{-2}}{F_n(\th)}\fibon{n}{k}\, e^{(n-1)\th}
[(1-q^k)(1-q^{n-k+1}) - q(1-q^{k-1})(1-q^{n-k})].
\end{gathered}
\nonumber
\eea
The quantity in the square bracket simplifies to $(1-q)(1-q^n)$
and the result follows.
\end{proof}

It is clear from \eqref{eqdefnFn} that
\bea
\label{eqincondfn}
F_1(\th) = 1, \quad F_2(\th) = 2\sinh \th.
\eea
We now choose $\th$ such that
\bea
\label{eqrestr}
 \sinh \th = \textup{a positive  integer}.
\eea It then follows from the three term recurrence relation for
Chebyshev polynomials that $\{F_n(\th)\}$ solves \eqref{eqfibrec}
under the initial conditions \eqref{eqincondfn}. This and
\eqref{eqrestr} show that $F_n(\th)$ is a positive integer for all
$n, n > 0$. Theorem 1.1 implies that $\fibon{n}{k}$ is always a
positive integer when $n >0$.

We can express \eqref{eqfib/fib} as the $n$th moment of the
measure \bea \nu(x) = (1-q^\al)\sum_{k=0}^\infty q^{\al k}
\delta(x- q^ke^{-\th}), \eea where $\delta(x-c)$ is a unit atomic
measure located at $x=c$. When $\al$ is even this is a positive
measure with total mass $=1$, otherwise $\nu$ is a unit signed
measure. In view of (18.4.11) and (18.4.13) in \cite{Ismbook} we
see that the corresponding orthogonal polynomials are little
$q$-Jacobi polynomials $\{p_n(xe^\theta; q^{\al-1}, 1) \}$,  where
\bea 
\begin{gathered}
p_n(x; a, b) = {}_2\phi_1(q^{-n}, abq^{n+1}; aq;q,qx) \\
= \sum_{j=0}^n \frac{(q;q)_n(abq^{n+1};q)_j}{(q;q)_j(q;q)_{n-j}}
q^{\binom{j+1}{2}} \frac{(-x)^j}{(aq;q)_j}, 
\end{gathered}
\eea
and the $q$-shifted factorials are 
$$(\la;q)_s = (1-\la)(1-\la q) \dots (1-\la q^{s-1}).$$
 The 
above $q$ is a base for the $q$-shifted factorials and is not 
the same as in \eqref{eqdfnq}

In terms of the generalized Fibonacci coefficients the polynomials
are expressed as
\bea
\label{eqdefpnal}
\begin{gathered}
p_n^{(\al)}(x):= \fibon{n+\al-1}{n} \;p_n(xe^\theta; q^{\al-1}, 1) \\
= \sum_{k=0}^n
\fibon{n}{k} \fibon{\al+n+k-1}{n}\, (-1)^{nk+\binom{k}{2}} \, x^k.
\end{gathered}
\eea The orthogonality relation is \cite{Gas:Rah}%
\bea
\label{eqortrelpnal}
\int_{\mathbb{R}}p_m^{(\al)}(x)p_n^{(\al)}(x)\, d\nu(x) =
(-1)^{\al n}\, \frac{F_\al(\th)}{F_{\al +n}(\th)}\, \delta_{m,n}.
\eea Recall that if the orthonormal polynomial of degree $n$ is
$$
\gamma_n x^n + \, \textup{lower order terms}
$$
then the Hankel determinant $D_n$ is given by
\bea D_n =
\prod_{j=1}^n \gamma_j^{-2}.
\eea
Consequently
\bea
\label{eqfibhd}
\begin{gathered}
\textup{det}\left\{1/F_{\al + i+j}(\th): 0 \le i,j \le n\right\}
 \qquad \qquad \\
= (-1)^{\al \binom{n+1}{2}} F_\al^{-n}(\th) \left[ \prod_{k=1}^n
F_{\al +2k}\fibon{\al+2k-1}{k} \right]^{-1}.
\end{gathered}
\eea
\begin{thm}
Let $A$ be the the matrix
$\{1/F_{\al + j+k}:0 \le j,k \le n\}$.
Then $A^{-1}$ has the matrix elements
\bea
\begin{gathered}
(-1)^{(\al+j+k)n- \binom{j}{2}-\binom{k}{2}} \;
F_{\al+j+k}(\th) \fibon{\al+n+j}{n-k} \\
\times \fibon{\al+n+k}{n-j} \; \fibon{\al+j+k-1}{j}\;
\fibon{\al+j+k-1}{k}.
\end{gathered}
\nonumber
\eea
\end{thm}
\begin{proof}
Use Lemma 1.1, \eqref{eqdefpnal}, and \eqref{eqortrelpnal}.
\end{proof}

\section{Generalized Lucas Numbers}
We now consider the Chebyshev polynomials of the first  kind
\bea
\begin{gathered}
T_n(i\sinh \th)= T_n(\cos(\pi/2-i\th) =
[e^{i(\pi/2-i\th)n} + e^{-i(\pi/2-i\th)n}]\\
= \frac{i^n}{2}\left [e^{n\th} + (-1)^ne^{-n\th}\right]
\end{gathered}
\nonumber \eea Define the generalized Lucas numbers by \bea
L_{n}(\th) = 2 (-i)^n \, T_n(\cos(i\th-\pi/2) =  \left[e^{n\th} +
(-1)^ne^{-n\th}\right], \eea for $n=0, 1, \dots$. Thus \bea
\label{eqincondln} L_0(\th) = 2, \quad L_1(\th) =  2\sinh \th.
\eea Hence $L_2(\th) = \cosh (2\th).$ It readily follows that
$\{L_n(\th)\}$ solves \eqref{eqfibrec} under the initial
conditions \eqref{eqincondln}. Assume \bea 2 \sinh \th = \textup{a
positive integer}. \eea Consequently $L_n(\theta)$ is a positive
integer for all $n$, $ n=1, 2, \dots$. Clearly there are
infinitely many such $\th$s. Moreover $L_n(\th_0) = L_n$.
% The
% notation for the Lucas numbers here is slightly different form
% the notation in \cite{Kos}. Koshy defines the Lucas sequence as
% $1, 3, 4, 7, 11, \dots$ while our sequence when $\th = \th_0$ is
% $2,1,3, 4, 7, \dots$, so the index is shifted by 1.

In view of  \eqref{eqdfnq} we see that %
\bea %
\label{eqlnasq}
L_{n}(\th) =  e^{n\th}\; [1+q^n]. %
\eea %
We extend the  definition of $L_n(\th)$ to $n \le 0$ by
\eqref{eqlnasq}. The explicit representation of $T_n(x)$, \cite[\S
4.5]{Ismbook} establishes the representation %
\bea%
 L_n(\th) = \sum_{k=0}^{\lfloor{n/2}\rfloor} \binom{n}{2k}
\sinh^{n-2k}( \th) \cosh^{2k}(\th). %
 \eea

It readily follows from \eqref{eqlnasq} that  %
\bea %
\frac{L_\al(\th)}{L_{n+\al}(\th)}
= (1+q^{\al})\sum_{k=0}^\infty (-q^{\al})^k(q^ke^{-\th})^n.
\nonumber%
 \eea %
 Define a measure $\psi$ by \bea \psi =
(1+q^\al)\sum_{k=0}^\infty (-q^{\al})^k \delta(x- q^ke^{-\th}),
\eea where, as before,  $\delta(x-c)$ is a unit atomic measure
located at $x=c$. Analogous to the definition
\eqref{eqdeffibbinom}
 the binomial coefficient relative to the generalized Lucas numbers
$\{L_n(\th)\}$ is
\bea
\label{eqdeflucbinom}
\luc{n}{0} := 1,
\quad \luc{n}{k}= \frac{L_n(\th) L_{n-1}(\th)\dots L_{n-k+1}(\th)}
{L_1(\th)\; L_2(\th)\;  \dots \; L_k(\th)}.
\eea
It is unlikely
that the binomial coefficients relative to the generalized Lucas
numbers are integers, but they may be integers if we only use the
generalized Lucas numbers of odd indices.

The polynomials
 \bea
\begin{gathered}
q_n^{(\al)}(x):= \luc{n+\al-1}{n} \;p_n(xe^\theta; -q^{\al-1}, 1) \\
= \sum_{k=0}^n \luc{n}{k} \luc{\al+n+k-1}{n}\,
(-1)^{nk+\binom{k}{2}} \, x^k
\end{gathered}
\eea
are special little $q$-Jacobi polynomials and satisfy the
orthogonality relation
\bea
\int_{\mathbb{R}}p_m^{(\al)}(x)p_n^{(\al)}(x)\, d\nu(x) =
(-1)^{\al n}\, \frac{L_\al(\th)}{L_{\al +n}(\th)}\, \delta_{m,n}.
\eea
 The proof of \eqref{eqfibhd} can be modified to establish
\bea
\begin{gathered}
\textup{det}\left\{1/L_{\al + i+j}(\th):0 \le i,j \le n\right\}  \qquad \qquad \\
= (-1)^{\al \binom{n+1}{2}} L_\al^{-n}(\th) \left[ \prod_{k=1}^n
L_{\al +2k}\luc{\al+2k-1}{k} \right]^{-1}.
\end{gathered}
\eea

\section{Relations}

Let $y_n$ be a solution to \eqref{eqfibrec} with integer
initial conditions and $y_1$ and $y_2$. If $y_0$ and $y_1$ are relatively
prime then $y_n$ and $y_{n+1}$ are relatively prime. This follows by
induction form \eqref{eqfibrec}. Consequently $F_n(\th)$ and
$F_{n+1}(\th)$ are relatively prime, and so are $L_n(\th)$ and
$L_{n+1}(\th)$.

The following result follows from \eqref{eqdfnq}--\eqref{eqfnasq}
and  \eqref{eqlnasq}.
\begin{thm}
For all  integers $\al, n, i,j$, the following identities  hold
 \bea
 \begin{gathered}
F_{\al + n +i}(\th)F_{\al + n +j}(\th) - (-1)^{\al + i + j} F_{ n
-i}(\th)F_{n -j}(\th) \\
= F_{\al +2n}(\th)F_{\al +i +j}(\th),
\end{gathered}
\eea %
\bea \label{equsedbyowi}%
F_{n+i}(\th) F_{n+j}(\th)-  F_{n}(\th)
F_{n+i+j}(\th) = (-1)^n F_i(\th)  F_{j}(\th),
 \eea
together with their companion formulas \bea
 \begin{gathered}
L_{\al + n +i}(\th)L_{\al + n +j}(\th) - (-1)^{\al + i + j}
(4+k^2) L_{n-i}(\th)L_{n -j}(\th) \\
= L_{\al +2n}(\th)L_{\al +i +j}(\th),
\end{gathered}
\eea %
\bea L_{n+i}(\th) L_{n+j}(\th)-  L_{n}(\th) L_{n+i+j}(\th) =
(-1)^{n+1} (k^2+4)F_i(\th)  F_{j}(\th),
 \eea
 where $k = 2 \sinh \th$.
\end{thm}

One interesting application of \eqref{equsedbyowi} is to take
$i=-j=2$, replace $n$ by $2n\pm 1$ and conclude that \bea
\label{eqowinscong}%
\qquad F_{2n+1}^2(\th) \equiv -k^2({\textup {mod}}\,
F_{2n-1}(\th)) \;\;  {\textup{and}} \;\; F_{2n-1}^2(\th) \equiv
-k^2 ({\textup {mod}}\, F_{2n+1}(\th)).
\eea
Thus given an integer $k$ a solution to the system of congruences
$$
a^2 \equiv -k^2 (\textup{mod} \; b), \quad \textup{and} \quad 
a^2 \equiv -k^2 (\textup{mod} \; b),
$$
is $(a,b) = (F_{2n-1} (\th), F_{2n+1} (\th))$. The converse to 
this may be   true, at least for certain values of  $k$ and it is 
interesting to characterize such values. In 
the  case $k=1$ the converse is due to Owings \cite{Owi}.

One of the topics in \S 32.3--32.4 in \cite{Kos} is the question
of evaluating the sums
$$
\sum_{i, j,k >0, i+j+k=n} F_iF_jF_k.
$$
We consider the more general question of evaluating $S_n$,
\bea
S_{m}(n) : = \sum_{j_1 , j_2, \dots, j_m: j_1 + j_2 + \dots + j_m =n} F_{j_1}F_{j_2}
\dots F_{j_m}.
\eea
Since $S_{m}(n)$ is an $m$-fold Cauchy convolution we find
$$
\Sum S_{m}(n) t^n = t^m(1-t-t^2)^{-m}.
$$
The ultraspherical polynomials polynomials $\{C_n^\nu(x)\}$ have the generating function
$$
\Sum C_n^\nu(x) t^n = (1-2xt+t^2)^{-\nu},
$$
\cite[\S4.5]{Ismbook}. They have the explicit formula
$$
C_n^\nu(x) = \sum_{k=0}^{\lfloor{n/2}\rfloor}
\frac{(2\nu)_n \, x^{n-2k}\, (x^2-1)^k}
{4^k\, k!\, (\nu +1/2)_k\, (n-2k)!}.
$$
It is clear that $S_{m}(n) =0$ if $n < m$. Therefore
\bea
\label{eqkillkosh}
\begin{gathered}
S_{m}(n+m) = (-i)^{n} C_n^m(i/2) \\
= \frac{(2m)_n}{2^n} \, \sum_{k=0}^{\lfloor{n/2}\rfloor}
\frac{(5/4)^k}{ k!\, (m+1/2)_k \, (n-2k)!}.
\end{gathered}
\eea 
Formula \eqref{eqkillkosh} generalizes many of the formulas
in \cite{Kos}. The only drawback of \eqref{eqkillkosh} is that it
does not show that $S_m(n)$ is an integer. Indeed the individual
terms in the sum (after multiplication by $2^{-n}(2m)_n$ are not
integers but their sum is an integer.

Theorem 5.9 in \cite{Kos} asserts that
$$
F_{n+k}F_{n-k} - F_n^2 = (-1)^{n+k+1}F_k^2,
$$
and is attributed to Catalan. Equations \eqref{eqdfnq} and
\eqref{eqfnasq} yield the identical result %
\bea
\label{eqgencassini} F_{n+k} (\th)\, F_{n-k}(\th) - F_n^2(\th) =
(-1)^{n+k+1}F_k^2(\th). %
\eea %
A consequence of \eqref{eqgencassini} is that if $p \mid F_n(\th)$
and $p \mid F_{n\pm k}(\th)$ then $p \mid F_k(\th)$. The case
$k=1$ of \eqref{eqgencassini} is
\bea \label{eqgencassini2} F_{n+1} (\th)\, F_{n-1}(\th) -
F_n^2(\th) = (-1)^{n}. \eea and generalizes the Cassini formula,
\cite[Theorem 5.3]{Kos}. Moreover %
\bea %
\label{eqmatrixfib}%
 \left(\begin{matrix} F_2(\th) \quad  F_1(\th)
\\
F_1(\th) \quad  F_0\th)  \end{matrix} \right)^n =
\left(\begin{matrix} F_{n+1}(\th) \quad F_n(\th)
\\
F_n(\th) \qquad F_{n-1}(\th)  \end{matrix}\right),
\eea
follows from equations \eqref{eqdfnq} and \eqref{eqfnasq}, and the
fact that%
\bea
 \left(\begin{matrix} F_2(\th) \quad
F_1(\th)
\\
F_1(\th) \quad F_0\th)  \end{matrix}\right) = \frac{1}{1+
e^{-2\th}} \left(\begin{matrix} 1  \qquad \;\; e^{-\th}
\\
e^{-\th} \quad  -1  \end{matrix}\right)  \left(\begin{matrix}
e^{\th}  \qquad  0
\\
0 \quad  -e^{-\th}  \end{matrix} \right)\left(\begin{matrix} 1
 \qquad  e^{-\th}
\\
e^{-\th} \quad   -1  \end{matrix} \right) \nonumber %
\eea %
One can prove a formula similar  to \eqref{eqmatrixfib} and
involving the generalized Lucas numbers. The relationship
\eqref{eqgencassini2} follows  from \eqref{eqmatrixfib} by
evaluating the determinants of both sides.

\begin{thm}
The generalized Fibonacci numbers have the property
\bea \label{eqinvtan} \qquad \arctan(k/F_{2m+1}(\th)) +
\arctan(1/F_{2m+2}(\th))= \arctan(1/F_{2m}(\th)),
\eea %
where $k = 2\sinh \th$. Moreover
\bea \label{eqsuminvtan}
\Sum \arctan(k/F_{2n+3}(\th)) = \arctan(1/k). %
\eea
\end{thm}
\begin{proof}
Clearly \eqref{eqinvtan} is equivalent to
$$ \left[\frac{k}{F_{2m+1}(\th)} + \frac{1}{F_{2m+2}(\th)}\right]/
\left[1 -\frac{k}{F_{2m+1}(\th)F_{2m+2}(\th)}\right] =
\frac{1}{F_{2m}(\th)}.
$$
In other words we need to show that
$$
[kF_{2m+2}(\th)+ F_{2m+1}(\th)] F_{2m}(\th)
 = F_{2m+1}(\th)F_{2m+2}(\th) -k.
 $$
 The above can be rewritten as
 $$
 k [ 1 + F_{2m+2}(\th) F_{2m}(\th)] = F_{2m+1}(\th) [F_{2m+2}(\th)- F_{2m}(\th)] =
 kF_{2m+1}^2(\th),
 $$
 which follows from \eqref{eqgencassini2}. Finally
 \eqref{eqsuminvtan} follows by telescopy from \eqref{eqinvtan}.
\end{proof}

It is easy to prove the following  result \bea \sum_{j=1}^n
F_j^2(\th) = 2 \sinh \th F_n(\th) F_{n+1}(\th). \eea which reduces
to a theorem of Lucas when $\th = \th_0$, see Theorem 5.5 in
\cite{Kos}. One can also prove \bea
\begin{gathered}
F_{n+1}^2(\th) + F_n^2(\th) = F_{2n+1}(\th), \\
F_{n+1}^2(\th) - F_n^2(\th) = 2 \sinh \th \, F_{2n+1}(\th).
\end{gathered}
\eea When $\th = \th_0$ the above identities reduce to results of
Lucas, \cite[Corollary 5.4]{Kos}.

The Lucas numbers are related to the Fibonacci numbers via
 \bea
\label{eqfn+fn=ln} L_{m}(\th) = F_{m+1}(\th)+ F_{m-1}(\th)
\eea
which follows from a calculation using \eqref{eqlnasq} and
\eqref{eqfnasq}.

The identities (5) and (7a) in \cite{Vaj} extend to
\bea
L_{n+1}(\th) + L_{n-1}(\th) = 4\cosh^2 \th F_n(\th)
\eea
and
\bea F_{n+2}(\th) - L_{n-2}(\th) = 2 \sinh \th L_n(\th), \eea
respectively.

Another identity which follows from \eqref{eqlnasq} and
\eqref{eqfnasq} is %
\bea %
\label{eqfmln}
 L_{m}(\th)F_{n}(\th) - F_{m+n}(\th) = (-1)^m
L_{n-m}(\th)
 \eea
 With $n = tm$ we iterate \eqref{eqfmln} and derive the finite
 continued fraction expansion
 \bea
\frac{F_{m(t+1)}(\th)}{F_{mt}(\th)} = L_m(\th) -
 \frac{(-1)^m}{L_m(\th)-}\; \frac{(-1)^m}{L_m(\th)-}\cdots
 \frac{(-1)^m}{L_m(\th)}
\eea %
In the above equation $L_m(\th)$ appears $m$ times.

\begin{thm}
The following identities hold %
\bea%
F_{2n}(\th) &=& F_n(\th)L_n(\th), \label{eqfnxln}\\
F_{n+m}(\th) &+& (-1)^m F_{n-m}(\th) = F_n(\th)L_m(\th)\label{eqfnxlm}\\
\sum_{j=0}^n \frac{1}{F_j(\th)}&=& \frac{1+\sinh \th}{\sinh \th} -
\frac{F_{2^n-1}(\th)}{F_{2^n}(\th)}. \label{eqsum1/f}
 \eea
\end{thm}
\begin{proof}
Formula \eqref{eqfnxln} is the special case $m=n$ of
\eqref{eqfnxlm}. The proof of \eqref{eqsum1/f} is by induction. It
clearly holds when $n=1$. The induction step uses
$$
-\frac{F_{2^n-1}(\th)}{F_{2^n}(\th)} + \frac{1}{F_{2^{n+1}}(\th)}
= -\frac{F_{2^n}(\th)}{F_{2^{n+1}}(\th)}
$$
which follows from \eqref{eqfnxln} and \eqref{eqfnxlm}.
\end{proof}

By letting $n \to \infty$ in Theorem 4.3 we find 
$$
\Sum \frac{1}{F_n(\th)} = 1 + e^{-\th}\coth \th.
$$
In the case of the Fibonacci numbers $\th = \th_0$ and the 
above sum reduces to (77) page 60 in \cite{Vaj}

\section{Integer Points on Algebraic Curves and Surfaces}
In this section we prove two theorems describing all the integral
points on the curves $y^2 -kxy -x^2 = \pm 1$ for a positive
integer $k$.
\begin{thm}
Let $\th > 0$ be given and assume that $k:=2sinh \th > 1$ is an
odd integer. A positive integer  $n$ is a generalized Fibonacci
number if and only if $4n^2\cosh^2 \th + 4$ or $4n^2\cosh ^2 \th
-4$ is a perfect square.
\end{thm}
\begin{proof} We will only consider
positive solutions to
\bea
x^2(k^2 + 4) - y^2 = \pm 4.
\label{eqgenpell}
\eea
It is clear that
\bea
\begin{gathered}
(k^2 + 4) F_m^2(\th) - L_m^2(\th) =
4\cosh ^2 \th F_m^2(\th) - L_m^2(\th) \\
= (e^\th + e^{-\th})^2
e^{2(m-1)\th}\left( \frac{1-q^m}{1-q}\right)^2- (1+q^m)^2e^{2m\th}
\end{gathered}
\nonumber
\eea
which simplifies to $4(-1)^{m+1}$, so it is equal to $\pm 1$. Hence
$n = F_m(\cosh \th)$ makes $4n^2\cosh^2 \th \pm 4$ a perfect square.
To prove the converse assume that
$4x_1^2 \cosh^2 \th \pm 4$ is
a perfect square $=y_1^2$ say. We assume $x_1 > 1$ and the case
$x_1 =1$ we considered at the end. Thus
$$y_1^2  = x_1^2(k^2+4) \pm 4  \quad \textup{with}\quad  x_1 > 1.
$$
It can be easily seen that $x_1>1$ implies $y_1 > kx_1$.
Let
$$
x_2 = (y_1-kx_1)/2, \quad y_2 = |(ky_1-(k^2+4)x_1)/2|.
$$
Both $x_2$ and $y_2$ are positive integers since $x_1$ and $y_1$ have
the same parity. A calculation shows that $x=x_2, y = y_2$ solve
\eqref{eqgenpell}. Moreover $x_2< x_1$ if and only if $y_1< (k+2)x_1$,
that is if and only if $x_1^2(k^2+4) \pm 4  < (k+2)^2 x_1^2$, since
the left-hand side is $y_1^2$. Clearly
the latter inequality holds, hence $0< x_2 < x_1$. We continue in this
manner until we reach $x_n =1$. Thus $y_n^2 = k^2 + 4 \pm 4$.The case $-$
leads to $y_n =k$ but the case $+$ makes $(y_n-k)(y_n+k) = 8$, hence
$y_n= k + 2^j$ and $y_n =-k + 2^{3-j}$, for some $j=0,1,2,3$. This forces
$y_n =(2^j +  2^{3-j})/2$ so that $j$ must equal 1 or 2, that
is $y_n =3$ which contradicts $k >1$. Thus the only solution is
$y_n=k = L_1(\th)$ and $x_1 = F_1(\th)$. By reversing the above steps,
and using  \eqref{eqfn+fn=ln} and \eqref{eqfibrec}
we see that $x_1 = F_n(\th)$ and $y_1 = L_n(\th)$.
\end{proof}
Note that in the process of proving Theorem 4.3 we also proved the
following.
\begin{cor}
We have
\bea
4\cosh ^2 F_n^2(\th) - L_n^2(\th) = 4(-1)^{n+1}
\eea
In particular $F_n(\th)\}$ and $\{L_n(\th)\}$ can not have any
common divisor larger than $2$. Moreover $F_n(\th)$ and $L_n(\th)$
have the same parity.
\end{cor}
Note that the Diophantine equation \eqref{eqgenpell} is a special
case of the Pell equation.

Let %
\bea
k := 2\sinh \th
\label{eqk-th}
 \eea
Observe that
\bea
\label{eq2ndalcu}
F_{n+1}^2(\th) - kF_n(\th)F_{n+1}(\th)-F_{n}^2(\th) = (-1)^n
\eea
follows from replacing the $F$'s in the left-hand side by the
corresponding expressions from \eqref{eqfnasq}. We now prove a
converse to \eqref{eq2ndalcu}. Consider the diophantine equations
\bea
&{}&y^2 -kxy - x^2 = 1, \label{eq+alcu} \\
&{}&y^2 -kxy - x^2 = -1.\label{eq-alcu}
\eea
The integer solutions to \eqref{eq+alcu} or \eqref{eq-alcu} will
be denoted by $(x,y)$. It is clear that if $(x,y)$ is such a pair
then $(-x,-y)$ will satisfy the same equation. Moreover if $(x,
y)$ satisfy \eqref{eq+alcu} (or \eqref{eq-alcu}) then $(y,-x)$
will solve \eqref{eq-alcu} (respectively \eqref{eq+alcu}). Hence
there is no loss of generality in assuming $x \ge 1$ and $y \ge
1$.
\begin{thm}
Let $k$ be an integer, $ k > 1$,  and related to $\th$ through
\eqref{eqk-th}. Assume that $(x,y)$ solve \eqref{eq+alcu}. Then
there exists a positive integer $n$ such that $(x, y) =
(F_{2n}(\th), F_{2n+1}(\th))$. On the other hand if $(x,y)$ solve
\eqref{eq-alcu} then there exists a positive integer $n$ such that
$(x, y) = (F_{2n-1}(\th), F_{2n}(\th))$.
\end{thm}
\begin{proof}The proof consists of three step.

\noindent{\bf Step 1}. We show that the smallest positive  $x$
satisfying \eqref{eq+alcu} is $x =k$. To see this write
\eqref{eq+alcu} in the form $y(y-kx) = x^2 +1$, hence $y = kx +z$
and $z >1$. Thus \eqref{eq+alcu} becomes $x(x-kz) = z^2 -1$, which
shows that $x \ge k = F_2(\th)$. The only possible answer for $y$
is $y = F_3(\th)$. Indeed the point $(F_2(\th), F_3(\th))$ lies on
the curve \eqref{eq+alcu}.

\noindent{\bf Step2} We use induction. Assume that all solutions
to \eqref{eq+alcu} are of the form $(F_{2j}(\th), F_{2j+1}(\th))$,
for $1 \le j \le m$. Let $x > F_{2m}(\th)$ and assume that $x$ is
the smallest integer such that $(x,y)$ solves \eqref{eq+alcu}.
Rewrite \eqref{eq+alcu} as
$$
(y-kx)^2 -1 = (k^2+1)x^2 -kxy = x[(k^2+1)x -ky].
$$
Thus $(k^2+1)x -ky > 0$. We have already shown that $y > kx$.
Define $(x_0,y_0)$ by
\bea
\label{eqdefx,y}
 x_0 = (k^2+1) x -ky,
\qquad y_0 = y -kx.
\eea
Both $x_0$ and $y_0$ are positive integers. Moreover $x_0 - x =
k(x-ky) <0$, that is $x_0 < x$, hence $x_0 \le F_{2m}(\th)$. By
direct computation we see that $(x_0, y_0)$ solves
\eqref{eq+alcu}, hence there is a positive integer $r$ such that
$x_0 = F_{2r}(\th)$ and $y_0 = F_{2r+1}(\th)$. From
\eqref{eqdefx,y} it follows that
$$
x = x_0+ ky_0, \quad \textup{and}\quad y = kx_0 + (1+k^2)y_0.
$$
Hence $x= F_{2r+2}(\th)$ and $y = F_{2r+3}(\th)$.

\noindent {\bf Step 3} Assume that $(x,y)$ solve \eqref{eq-alcu}
and set $(x_0, y_0) = (y, x+ky)$. A calculation shows that $(x_0,
y_0)$ satisfies \eqref{eq+alcu}, hence $(y, x+ky) =  (F_{2j}(\th),
F_{2j-1}(\th))$, for some positive integer $j$, which implies
$(x,y) = (F_{2j}(\th), F_{2j-1}(\th))$,  and the proof is
complete.
\end{proof}

We next extend the following identities of Carlitz \cite[Ex
91-91]{Kos}:
\bea%
\label{eqCar}
\begin{gathered}
Z_{n+1}^3- Z_n^3 - Z_{n-1}^3 = 3 Z_{n+1}Z_nZ_{n-1}, \quad Z_j =
F_j \; \textup{or}\;  L_j.
\end{gathered}
\nonumber \eea
\begin{thm}
With $k = \sinh \th$ the identity
 \bea
 \label{eqgeCar}
\begin{gathered}
Z_{n+1}^3(\th)-k^3 Z_n^3(\th) - Z_{n-1}^3 (\th) \\
=  3k Z_{n+1}(\th)Z_n (\th)Z_{n-1}(\th),
\end{gathered}
\eea %
 holds for $Z_n(\th) = F_n(\th)$ or $Z_n(\th) = L_n(\th)$.
\end{thm}
\begin{proof}After using \eqref{eqfibrec} we see that
the left-hand side of the above equation in the Fibonacci case is
\bea
\begin{gathered}
2\sinh \th F_n(\th) [F_{n+1}^2(\th)+ F_{n-1}^2 (\th)
+F_{n+1}(\th)F_{n-1}(\th)]
-(2\sinh \th)^3 F_n^3(\th) \\
= 2\sinh \th F_n(\th)  F_{n-1} (\th) [F_{n+1}(\th)+ 2\sinh \th
F_n(\th) + F_{n-1} (\th) + F_{n+1}(\th)]
\end{gathered}
\nonumber \eea which simplifies to the right-hand side of
\eqref{eqgeCar}. We only used the recurrence  relation
\eqref{eqfibrec} to establish \eqref{eqgeCar}. Thus
\eqref{eqgeCar} also holds for $\{L_n(\th)\}$ since it also
satisfies \eqref{eqfibrec}.
\end{proof}

It is interesting to determine all the positive integer points on
the surface  $z^3 - y^3 - z^3 = 3xyz$. We suspect that the only
solutions are $(x, y, z) = (F_{n-1}, F_n, F_{n+1})$ or $(L_{n-1},
L_n, L_{n+1})$. This would give a converse to Carlitz's identities
\eqref{eqCar}. Similarly it is of interest to determine all the
the positive integer points $(x, y, z)$ which lie on the surface
$z^3 - k^3y^3 - z^3 = 3k xyz$ for a given positive integer $k$.

Fairgrieve and Gould \cite{Fai:Gou} studied formulas involving
differences of products of Fibonacci numbers. They claim that
computer searches yielded only the list of formulas stated below.
They pointed out that
some of these formulas were already known and references are given in \cite{Fai:Gou}. %
\bea%
\label{eqfair-goul1}
F_{n+1} F_{n+2} F_{n+6}- F_{n+3}^2 &=& (-1)^n F_{n},
\\
F_{n} F_{n+4} F_{n+5} - F_{n+1}^3 &=& (-1)^{n+1} F_{n+6},\\
F_{n-2}F_{n+1}^2- F_n^3 &=& (-1)^{n-1}  F_{n-1}, %
\\
F_{n+2}F_{n-1}^2 - F_n^3 &=& (-1)^n F_{n+1}, \\
F_{n-3}F_{n+1}^3 - F_n^4 &=& (-1)^n \left[F_{n-1}F_{n+3}+
2F_n^2\right] \\
 F_{n+3}F_{n-1}^3- F_{n}^4 &=& (-1)^n \left[F_n^2+ F_n F_{n-1} + 2
 F_{n-1}^2\right].
\eea%
It is clear that we can rewrite the last equation above as %
\bea
F_{n+3}F_{n-1}^3- F_{n}^4 &=& (-1)^n \left[F_n F_{n+1}+  2
 F_{n-1}^2\right]
\eea
 These
can be extended to the numbers $F_n(\th)$. The
extensions are given below. %
\bea %
&{}& F_{n+1}(\th) F_{n+2}(\th) F_{n+6}(\th) - F_{n+3}^2(\th) \\
 &{}& \qquad =
(-1)^n \left[k^2 F_{n}(\th) + (k^3-1) F_{n+1}(\th)\right],
\nonumber \\
&{}& F_{n}(\th) F_{n+4}(\th) F_{n+5}(\th)
- F_{n+1}^3(\th) \\
&{}& \qquad = (-1)^{n+1} \left[F_{n+6}(\th) + k(k-1) F_{n+4}(\th)\right], \nonumber
\\
 &{}& F_{n-2}(\th)F_{n+1}^2(\th)- F_n^3 (\th)= (-1)^{n-1} \left[k F_{n-1}(\th) + (k^2-1) F_n(\th)\right], \\
 &{}& F_{n+2}(\th)F_{n-1}^2(\th) - F_n^3(\th) = (-1)^n [F_{n}(\th) + k F_{n-1}(\th)], \\
 &{}& F_{n-3}(\th)F_{n+1}^3(\th) - F_n^4(\th) \\
 &{}&= (-1)^n \left[F_{n-1}(\th)F_{n+3}(\th)+
2F_n^2(\th)+(k^2-1) F_n(\th)F_{n+2}(\th)\right] \nonumber \\
 &{}&  F_{n+3}(\th)F_{n-1}^3(\th)- F_{n}^4(\th) \\
 &{}& = (-1)^n \left[F_n^2(\th)+ F_n(\th) F_{n-1}(\th) + 2
 F_{n-1}^2(\th)\right].\nonumber
\eea %
The proofs use \eqref{eqgencassini}--\eqref{eqgencassini2} and
\eqref{eqfibrec}.

 \bigskip

{\bf Acknowledgements}: I wish to thank my friends Richard Askey, 
Christian Berg, abd Edwin Clark. Askey 
pointed out reference \cite{Vaj}. Berg  provided  me with a 
copy of his work  \cite{Ber} which
initiated my interest in the subject. Edwin Clark send me a copy 
of \cite{Kal:Men} and made very interesting remarks. 
This work was done while the
author was visiting the Liu Bie Ju center for Mathematical
Sciences of the City University of Hong Kong and he acknowledges
the hospitality and financial support.

% \bigskip

\end{document}